\documentclass[11pt,reqno]{amsart}
\usepackage{enumerate}
\usepackage{amsrefs}
\usepackage{amsfonts, amsmath, amssymb, amscd, amsthm, bm, cancel,mathrsfs}
\usepackage{url}
\usepackage{graphicx}
\usepackage[%backref=page,
linktocpage=true,colorlinks,citecolor=magenta,linkcolor=blue,urlcolor=magenta]{hyperref}
\usepackage{multicol}
\usepackage{comment}
\usepackage[margin=1in]{geometry}
\usepackage{tikz,pgfplots}

\pgfplotsset{compat=1.16}
\usepgfplotslibrary{fillbetween}

\newtheorem{thm}{Theorem}[section]
\newtheorem{cor}[thm]{Corollary}

\newtheorem*{conj}{Conjecture}

\newtheorem{prop}[thm]{Proposition}
\theoremstyle{definition}

\theoremstyle{remark}
\newtheorem{rmk}[thm]{Remark}

\numberwithin{equation}{section}
\allowdisplaybreaks

\newcommand{\bR}{\mathbb{R}}

\newcommand{\la}{\langle}
\newcommand{\ra}{\rangle}

\newcommand{\divg}{\mathrm{div}}

\begin{document}
\title[uniqueness of positive harmonic functions]{A proof of Guo-Wang's conjecture on the uniqueness of positive harmonic functions in the unit ball}
	
\author{Pingxin Gu}
	\address{Department of Mathematical Sciences, Tsinghua University, Beijing 100084, P.R. China}
	\email{\href{mailto:gpx21@mails.tsinghua.edu.cn}{gpx21@mails.tsinghua.edu.cn}}
\author{Haizhong Li}
	\address{Department of Mathematical Sciences, Tsinghua University, Beijing 100084, P.R. China}
	\email{\href{mailto:lihz@tsinghua.edu.cn}{lihz@tsinghua.edu.cn}}
	
%	\date{\today}

\begin{abstract}
	Guo-Wang [{\em Calc.Var.Partial Differential Equations}, {\bf 59} (2020)] conjectured that for $1<q<\frac{n}{n-2}$ and $0<\lambda\leq \frac{1}{q-1}$, the positive solution $u\in C^{\infty}(\bar B)$ to the equation
	\[
		\left\{
		\begin{array}{ll}
		\Delta u=0 &in\ B^n,\\
		u_{\nu}+\lambda u=u^q&on\ S^{n-1},
		\end{array}
	\right.
	\]
	must be constant. In this paper, we give a proof of this conjecture.
\end{abstract}

\maketitle
%	\tableofcontents

\section{Introduction}
In the past decades, a great deal of mathematical effort in best constant of Sobolev inequality has been devoted. For $n\geq 3$, a well-known subject is to figure out the minimum constant of Sobolev trace inequalities 
\begin{align}
||u||_{L^{\frac{2(n-1)}{n-2}}(\partial\bR_+^n)}\leq C||\nabla u||_{L^2(\bR_+^n)},\quad \forall u\in C_0^{\infty}(\bar \bR_+^n).
\end{align}
A key issue for this study is to investigate the extreme value of Sobolev quotient. Escobar {\cite{Esc88} showed by conformal transformations that the best constant equals to
\begin{align}\label{Sobolev quotient}
Q(B^n):=\inf_{u\in C^{\infty}(\bar B^n)}\frac{\int_{B^n}|\nabla u|^2+\frac{n-2}{2}\int_{S^{n-1}}u^2}{\big(\int_{S^{n-1}}|u|^{\frac{2(n-1)}{n-2}}\big)^{\frac{n-2}{n-1}}}.
\end{align}
Lions \cite{Lio85} proved that \eqref{Sobolev quotient} can be achieved by a positive $u$ satisfying the Euler-Lagrange equation
\begin{align}\label{E-L equation}
\left\{\begin{array}{ll}
\Delta u=0 &in\ B^n,\\
u_{\nu}+\frac{n-2}{2}u=u^{\frac{n}{n-2}}&on\ S^{n-1},
\end{array}\right.
\end{align}
where $\nu$ is the unit outer normal vector on $S^{n-1}$. With this conclusion, Escobar \cite{Esc88} classified all positive solutions of \eqref{E-L equation} by an integral method and hence \cite{Esc90} proved that
\begin{align}\label{Escobar inequality}
|S^{n-1}|^{\frac{1}{n-1}}\bigg(\int_{S^{n-1}}u^{\frac{2(n-1)}{n-2}}\bigg)^{\frac{n-2}{n-1}}\leq \frac{2}{n-2}\int_{B^n}|\nabla u|^2+\int_{S^{n-1}}u^2,\quad \forall u\in C^{\infty}(\bar B^n).
\end{align}
Different from Escobar, using harmonic analysis, Beckner \cite{Bec93} derived a family of inequalities
\begin{align}\label{Beckner inequality}
|S^{n-1}|^{\frac{q-1}{q+1}}\bigg(\int_{S^{n-1}}u^{q+1}\bigg)^{\frac{2}{q+1}}\leq (q-1)\int_{B^n}|\nabla u|^2+\int_{S^{n-1}}u^2,\quad \forall u\in C^{\infty}(\bar B^n),
\end{align}
provided $1<q<\infty$, if $n=2$, and $1<q\leq \frac{n}{n-2}$, if $n\geq 3$. The corresponding Euler-Lagrange equation to \eqref{Beckner inequality} is
\begin{align}\label{E-L equation-2}
\left\{\begin{array}{ll}
\Delta u=0 &in\ B^n,\\
u_{\nu}+\frac{1}{q-1}u=u^q&on\ S^{n-1}.
\end{array}\right.
\end{align}
It is apparent that the case $n\geq 3$ and $q=\frac{n}{n-2}$ of \eqref{Beckner inequality} and \eqref{E-L equation-2} are just \eqref{Escobar inequality} and \eqref{E-L equation} respectively. Also, in the same paper, Beckner \cite{Bec93} confirmed
\begin{align}\label{Beckner ineq}
|S^{n-1}|^{\frac{q-1}{q+1}}\bigg(\int_{S^{n-1}}u^{q+1}\bigg)^{\frac{2}{q+1}}\leq \frac{q-1}{n-1}\int_{S^{n-1}}|\nabla u|^2+\int_{S^{n-1}}u^2,\quad \forall u\in C^{\infty}(S^{n-1}),
\end{align}
provided $1<q<\infty$, if $n=2$ or $3$, and $1<q\leq \frac{n+1}{n-3}$, if $n\geq 4$. By considering the Euler-Language equation and using integral method, Bidaut-V\'eron and V\'eron \cite{BM91} gave a new proof of \eqref{Beckner ineq}. \\
\indent A natural question is: Now that \eqref{Beckner ineq} can be proved by the method of integration, can one prove \eqref{Beckner inequality} with the same strategy? Inspired by the arguments, Guo-Wang \cite{GW20} proposed the following conjecture.
\begin{conj}[\cite{GW20}]\label{conj}
If $u\in C^{\infty}(\bar B^n)$ is positive solution of the following equation
\[
	\left\{
		\begin{array}{ll}
		\Delta u=0 &in\ B^n,\\
		u_{\nu}+\lambda u=u^q&on\ S^{n-1}.
		\end{array}
	\right.
\]
Then $u$ is constant provided $1<q<\frac{n}{n-2}$ and $0<\lambda\leq \frac{1}{q-1}$.
\end{conj}
In recent years, there are some partial results about the conjecture, see \cite{GHW21,GW20,LO23}. A remarkable one is that Guo-Hang-Wang \cite{GHW21} confirmed the conjecture in $n=2$. In this paper, we investigate satisfactory Obata type identities. Combine with auxiliary functions, we give a proof of the conjecture when $n\geq 3$. 
\begin{thm}\label{main thm}
For $n\geq 3$, suppose $u\in C^{\infty}(\bar B^n)$ is a positive solution of the following equation
\begin{align}\label{initial equation}
	\left\{
		\begin{array}{ll}
		\Delta u=0 &in\ B^n,\\
		u_{\nu}+\lambda u=u^q&on\ S^{n-1}.
		\end{array}
	\right.
\end{align}
If $1<q\leq \frac{n}{n-2}$ and $0<\lambda\leq \frac{1}{q-1}$, then $u$ must be constant $\lambda^{\frac{1}{q-1}}$ unless $q=\frac{n}{n-2}$ and $\lambda=\frac{1}{q-1}$, $u$ is given by the following formula
\begin{align}
u_{\xi}(x)=\bigg(\frac{n-2}{2}\frac{1-|\xi|^2}{1+|\xi|^2|x|^2-2\la \xi,x\ra}\bigg)^{\frac{n-2}{2}},
\end{align}
for some $\xi\in B^n$.
\end{thm}

From Theorem \ref{main thm}, we immediately get
\begin{cor}
The conjecture holds for $n\geq 3$.
\end{cor}

\indent The paper is organized as follows. In Sect.\ref{sec:2}, we establish a Pohozaev identity by calculating divergence of a given tensor field with two parameters, and designate one of the parameters to make Pohozaev identity works. In Sect.\ref{sec:3}, we give an Obata type identity and determine the other parameter to ensure the effectiveness of Obata's skill. In Sect.\ref{sec:4}, we introduce auxiliary functions involving the length of position vector with one parameter to obtain improved identities. In Sect.\ref{sec:5}, we give the proof of Theorem \ref{main thm} and a new proof of Beckner's inequality \eqref{Beckner inequality} by integral method.

\section{A Pohozaev identity}\label{sec:2}
In this section, we use divergence theorem to prove an identity, and use Pohozaev identity to simplify it. During the simplification, we will determine the parameter as what we have hoped.
\subsection{Preliminary}
To begin with, let $u=v^{-\frac{1}{q-1}}$ in \eqref{initial equation}, then $v$ satisfies
\begin{align}\label{main equation}
	\left\{
		\begin{array}{ll}
		\Delta v=\frac{q}{q-1}\frac{|\nabla v|^2}{v}&in\ B^n,\\
		v_{\nu}=(q-1)(\lambda v-1)&on\ S^{n-1}.
		\end{array}
	\right.
\end{align}
\indent In this case, the boundary condition is concise, which only involves normal derivatives and linear expression of the solution $v$. In the text that follows, we define
\begin{align}\label{def of MN}
M:=\divg(v^a\nabla_{\nabla v}\nabla v),\qquad N:=\divg(v^a\Delta v \nabla v),
\end{align}
and focus on the quantity $M-bN$, where parameters $a,b\in \bR$ are to be determined later. \\
\indent In order to integrate by parts and handle the boundary term, we choose an orthonormal frame $\{e_i\}_{i=1}^n$ such that $e_n=\nu$ is the unit outward normal vector on $S^{n-1}$ and the second fundamental form of $S^{n-1}$ equals to identity. Thus on $S^{n-1}$, by use of Reilly's formula, see \cite{Rei77}, we have
\begin{align}
\label{boundary-1}	v_{\alpha n}=&\big(\lambda(q-1)-1\big)v_{\alpha},\qquad \forall\ 1\leq \alpha\leq n-1,\\
\label{boundary-2}	\sum_{\alpha=1}^{n-1}v_{\alpha\alpha}=&\Delta_{S^{n-1}}v+(n-1)v_n,
\end{align}
and
\begin{align}
\label{boundary-3}	v_{nn}=&\Delta v-\sum_{\alpha=1}^{n-1}v_{\alpha\alpha}=\frac{q}{q-1}\frac{|\nabla v|^2}{v}-\Delta_{S^{n-1}}v-(n-1)v_n.
\end{align}
For simplicity of presentation, we always set subscripts $1\leq \alpha\leq n-1$ and $1\leq i,j\leq n$. Einstein summation convention for these subscripts is always used in what follows.
\subsection{Choice of parameter a}
With these preliminaries, we integrate $M$ and $N$ in $B^n$ respectively. Combining with \eqref{boundary-1}, \eqref{boundary-3} and \eqref{main equation}, we obtain by use of divergence theorem
\begin{align}\label{Integration on M}
	\int_{B^n}M=&\int_{B^n}(v^av_{ij}v_i)_j=\int_{S^{n-1}}v^av_{nn}v_n+\int_{S^{n-1}}v^av_{\alpha n}v_{\alpha}\nonumber\\
	=&\frac{q}{q-1}\int_{S^{n-1}}v^{a-1}|\nabla v|^2v_n-\int_{S^{n-1}}v^a\Delta_{S^{n-1}}vv_n-(n-1)\int_{S^{n-1}}v^av_n^2\nonumber\\
	&+\Big(\lambda(q-1)-1\Big)\int_{S^{n-1}}v^a|\nabla_{S^{n-1}} v|^2,\\
\label{Integration on N}	\int_{B^n}N=&\int_{B^n}(v^a\Delta vv_i)_i=\int_{S^{n-1}}v^a\Delta vv_n=\frac{q}{q-1}\int_{S^{n-1}}v^{a-1}|\nabla v|^2v_n.
\end{align}
Now we deal the term $\int_{S^{n-1}}v^a\Delta_{S^{n-1}}vv_n$ in \eqref{Integration on M} by applying divergence theorem on $S^{n-1}$:
\begin{align}\label{Divergence thm on M}
	\int_{S^{n-1}}v^a\Delta_{S^{n-1}}vv_n=&(q-1)\lambda\int_{S^{n-1}}v^{a+1}\Delta_{S^{n-1}}v-(q-1)\int_{S^{n-1}}v^a\Delta_{S^{n-1}}v\nonumber\\
	=&-(q-1)(a+1)\lambda\int_{S^{n-1}}v^a|\nabla_{S^{n-1}}v|^2+(q-1)a\int_{S^{n-1}}v^{a-1}|\nabla_{S_{n-1}}v|^2.
\end{align}
Combining formulas \eqref{Integration on M}, \eqref{Integration on N} through \eqref{Divergence thm on M} and splitting $|\nabla v|^2$ into $|\nabla_{S^{n-1}}v|^2+v_n^2$ on $S^{n-1}$, we derive
\begin{align}\label{Integration on Q2}
\begin{aligned}
	\int_{B^n}(M-bN)=&(1-b)\frac{q}{q-1}\int_{S^{n-1}}v^{a-1}|\nabla_{S^{n-1}}v|^2v_n+(1-b)\frac{q}{q-1}\int_{S^{n-1}}v^{a-1}v_n^3\\
	&+\lambda(q-1)(a+1)\int_{S^{n-1}}v^a|\nabla_{S^{n-1}}v|^2-(q-1)a\int_{S^{n-1}}v^{a-1}|\nabla_{S^{n-1}}v|^2\\
	&-(n-1)\int_{S^{n-1}}v^av_n^2+\Big(\lambda(q-1)-1\Big)\int_{S^{n-1}}v^a|\nabla_{S^{n-1}}v|^2.
\end{aligned}
\end{align}
Use \eqref{main equation} to eliminate $v_n$ in the term $(1-b)\frac{q}{q-1}\int_{S^{n-1}}v^{a-1}|\nabla_{S^{n-1}}v|^2v_n$ and one of the three $v_n$'s in the term $(1-b)\frac{q}{q-1}\int_{S^{n-1}}v^{a-1}v_n^3$, \eqref{Integration on Q2} becomes
\begin{align}\label{Integration on Q3}
\begin{aligned}
	\int_{B^n}(M-bN)=&\lambda q(1-b)\int_{S^{n-1}}v^a|\nabla_{S^{n-1}}v|^2-q(1-b)\int_{S^{n-1}}v^{a-1}|\nabla_{S^{n-1}}v|^2\\
	&+\lambda q(1-b)\int_{S^{n-1}}v^av_n^2-q(1-b)\int_{S^{n-1}}v^{a-1}v_n^2\\
	&+\lambda(q-1)(a+1)\int_{S^{n-1}}v^a|\nabla_{S^{n-1}}v|^2-(q-1)a\int_{S^{n-1}}v^{a-1}|\nabla_{S^{n-1}}v|^2\\
	&-(n-1)\int_{S^{n-1}}v^av_n^2+\Big(\lambda(q-1)-1\Big)\int_{S^{n-1}}v^a|\nabla_{S^{n-1}}v|^2.
\end{aligned}
\end{align}
We obtain from \eqref{Integration on Q3} that
\begin{align}\label{Integration on Q}
\begin{aligned}
	\int_{B^n}(M-bN)=&\Big(\lambda q(1-b)+\lambda(q-1)(a+2)-1\Big)\int_{S^{n-1}}v^a|\nabla_{S^{n-1}}v|^2\\
	&-\Big(q(1-b)+a(q-1)\Big)\int_{S^{n-1}}v^{a-1}|\nabla_{S^{n-1}}v|^2\\
	&+\Big(\lambda q(1-b)-(n-1)\Big)\int_{S^{n-1}}v^av_n^2\\
	&-q(1-b)\int_{S^{n-1}}v^{a-1}v_n^2.
\end{aligned}
\end{align}
The term $\int_{S^{n-1}}v^{a-1}|\nabla_{S^{n-1}}v|^2$ in \eqref{Integration on Q} appears for the reason that the boundary condition of \eqref{main equation} is not homogeneous. It is desirable to eliminate it with some equalities. The key idea is to confirm Pohozaev identities from conditions \eqref{main equation}. Only if we choose $a=-\frac{q+1}{q-1}$ makes it work.

\begin{prop}\label{Pohozaev identity}
Let $v$ be a positive solution of \eqref{main equation}. For $a=-\frac{q+1}{q-1}$, we derive the following Pohozaev identity
\begin{align}\label{Pohozaev-0}
\begin{aligned}
\int_{S^{n-1}}&v^{a-1}|\nabla_{S^{n-1}}v|^2=\int_{S^{n-1}}v^{a-1}v_n^2+(n-2)\int_{S^{n-1}}v^av_n^2-(n-2)(q-2)\lambda\int_{B^n}v^a|\nabla v|^2.
\end{aligned}
\end{align}
\end{prop}
\noindent{\it Proof.} We note $x_{ij}=\delta_{ij}$ in $B^n$. Thus, we have the following calculation
\begin{align}\label{Pohozaev-1}
\begin{aligned}
	\int_{B^n}\divg(v^{a-1}|\nabla v|^2&x)=\\
	&(a-1)\int_{B^n}v^{a-2}|\nabla v|^2x_iv_i+2\int_{B^n}v^{a-1}v_{ij}x_iv_j+n\int_{B^n}v^{a-1}|\nabla v|^2.
\end{aligned}
\end{align}
On the other hand, take \eqref{main equation} into consideration, we obtain
\begin{align}\label{Pohozaev-1.5}
\begin{aligned}
	\int_{B^n}\divg(v^{a-1}&\la \nabla v,x\ra \nabla v)=\\
	&\Big(a-1+\frac{q}{q-1}\Big)\int_{B^n}v^{a-2}|\nabla v|^2x_iv_i+\int_{B^n}v^{a-1}v_{ij}x_iv_j+\int_{B^n}v^{a-1}|\nabla v|^2.
\end{aligned}
\end{align}
Note the condition $a=-\frac{q+1}{q-1}$ is the only choice to satisfy
\[
\frac{a-1}{a-1+\frac{q}{q-1}}=2.
\]
Then \eqref{Pohozaev-1}$-2\times$\eqref{Pohozaev-1.5} implies
\begin{align}\label{Pohozaev-1.7}
(n-2)\int_{B^n}v^{a-1}|\nabla v|^2=\int_{B^n}\divg(v^{a-1}|\nabla v|^2x)-2\int_{B^n}\divg(v^{a-1}\la \nabla v,x\ra \nabla v).
\end{align}
Notice on $S^{n-1}$, we have $x=\nu$, thus $\la x,\nu\ra=1$ and $\la \nabla v,x\ra=\la \nabla v,\nu\ra=v_n$. By divergence theorem, we have from \eqref{Pohozaev-1.7} that
\begin{align}\label{Pohozaev-2}
\begin{aligned}
(n-2)\int_{B^n}v^{a-1}|\nabla v|^2=&\int_{S^{n-1}}v^{a-1}|\nabla v|^2\la x,\nu\ra-2\int_{S^{n-1}}v^{a-1}\la \nabla v,x\ra\la \nabla v,\nu\ra\\
=&\int_{S^{n-1}}v^{a-1}|\nabla v|^2-2\int_{S^{n-1}}v^{a-1}v_n^2\\
=&\int_{S^{n-1}}v^{a-1}|\nabla_{S^{n-1}} v|^2-\int_{S^{n-1}}v^{a-1}v_n^2.
\end{aligned}
\end{align}
Now it remains to deal with the term $\int_{B^n}v^{a-1}|\nabla v|^2$ on the left-hand side of \eqref{Pohozaev-2}. A key observation follows from divergence theorem and $a=-\frac{q+1}{q-1}$ that
\begin{align}\label{Pohozaev-3}
\int_{S^{n-1}}v^av_n=\int_{B^n}\divg(v^a\nabla v)=\Big(a+\frac{q}{q-1}\Big)\int_{B^n}v^{a-1}|\nabla v|^2=-\frac{1}{q-1}\int_{B^n}v^{a-1}|\nabla v|^2,
\end{align}
and
\begin{align}\label{Pohozaev-3.5}
\int_{S^{n-1}}v^{a+1}v_n=\int_{B^n}\divg(v^{a+1}\nabla v)=\Big(a+1+\frac{q}{q-1}\Big)\int_{B^n}v^a|\nabla v|^2=\frac{q-2}{q-1}\int_{B^n}v^a|\nabla v|^2.
\end{align}
Using the fact that the boundary condition of \eqref{main equation} yields
\begin{align}\label{Pohozaev-4}
\int_{S^{n-1}}v^av_n^2=(q-1)\int_{S^{n-1}}v^av_n(\lambda v-1)=\lambda(q-1)\int_{S^{n-1}}v^{a+1}v_n-(q-1)\int_{S^{n-1}}v^av_n.
\end{align}
Putting \eqref{Pohozaev-3} and \eqref{Pohozaev-3.5} into \eqref{Pohozaev-4}, we conclude
\begin{align}\label{Pohozaev-5}
\int_{S^{n-1}}v^av_n^2=\lambda(q-2)\int_{B^n}v^a|\nabla v|^2+\int_{B^n}v^{a-1}|\nabla v|^2.
\end{align}
Then \eqref{Pohozaev-0} follows from substituting \eqref{Pohozaev-5} into \eqref{Pohozaev-2} to eliminate the term $\int_{B^n}v^{a-1}|\nabla v|^2$.
\\\qed\\
\indent By means of Proposition \ref{Pohozaev identity}, we eliminate $\int_{S^{n-1}}v^{a-1}|\nabla_{S^{n-1}}v|^2$ in \eqref{Integration on Q} and derive
\begin{cor}
Let $v$ be a positive solution of \eqref{main equation}. The quantities $M,N$ are defined as \eqref{def of MN}. Then for $a=-\frac{q+1}{q-1}$, 
\begin{align}\label{Integration on Q-2}
\begin{aligned}
	\int_{B^n}(M-bN)=&\Big(\lambda q(1-b)+\lambda(q-3)-1\Big)\int_{S^{n-1}}v^a|\nabla_{S^{n-1}}v|^2\\
	&+\Big(q(1-b)-(q+1)\Big)(n-2)(q-2)\lambda\int_{B^n}v^a|\nabla v|^2\\
	&+\Big(\lambda q(1-b)+(n-2)qb-1\Big)\int_{S^{n-1}}v^av_n^2\\
	&-\Big(2q(1-b)-(q+1)\Big)\int_{S^{n-1}}v^{a-1}v_n^2.
\end{aligned}
\end{align}
\end{cor}  

In the next section, we will determine the value of $b$, which is a rather subtle issue. One can guess that we should control the coefficient of $\int_{S^{n-1}}v^{a-1}v_n^2$ as
\begin{align}\label{c geq 0}
2q(1-b)-(q+1)\geq 0,
\end{align}
since we have no other ways to decompose it. In fact, we will make the equality sign in \eqref{c geq 0} hold. More logical reasons for such choice will be explained later. 

\section{On Obata's skill}\label{sec:3}
\subsection{Obata type identities}
It is time to expand $M-bN$, which gives
\begin{align}\label{M-bN ob1}
\begin{aligned}
	M-bN=&(v^av_{ij}v_i)_j-b(v^a\Delta vv_j)_j\\
	=&\Big(v^av_{ij}v_{ij}+av^{a-1}v_{ij}v_iv_j+v^a(\Delta v)_iv_i\Big)-b\Big(av^{a-1}\Delta v|\nabla v|^2+v^a(\Delta v)_jv_j+v^a(\Delta v)^2\Big)\\
	=&v^av_{ij}v_{ij}+av^{a-1}v_{ij}v_iv_j+(1-b)v^a(\Delta v)_iv_i-abv^{a-1}\Delta v|\nabla v|^2-bv^a(\Delta v)^2.
	\end{aligned}
\end{align}
Recall that $\Delta v$ satisfies \eqref{main equation}, so
\begin{align}\label{M-bN ob2}
v^a(\Delta v)_iv_i=\frac{q}{q-1}v^a\Big(\frac{|\nabla v|^2}{v}\Big)_iv_i=\frac{2q}{q-1}v^{a-1}v_{ij}v_iv_j-\frac{q}{q-1}v^{a-2}|\nabla v|^4.
\end{align}
If we put \eqref{M-bN ob2} and \eqref{main equation} into \eqref{M-bN ob1}, we will obtain
\begin{align}\label{M-bN ob3}
M-bN=&v^av_{ij}v_{ij}+\Big(a+(1-b)\frac{2q}{q-1}\Big)v^{a-1}v_{ij}v_iv_j\nonumber\\
&-\Big((1-b)\frac{q}{q-1}+ab\frac{q}{q-1}+b(\frac{q}{q-1})^2\Big)v^{a-2}|\nabla v|^4.
\end{align}
For simplicity, we define $d$ as
\[
d:=\frac{a}{2}+(1-b)\frac{q}{q-1}=\frac{2q(1-b)-(q+1)}{2(q-1)}.
\]
After putting $a=-\frac{q+1}{q-1}$, \eqref{M-bN ob3} becomes
\begin{align}\label{M-bN ob4}
M-bN=&v^av_{ij}v_{ij}+2dv^{a-1}v_{ij}v_iv_j-\Big(\frac{q}{q-1}d+\frac{q}{2(q-1)}\Big)v^{a-2}|\nabla v|^4\nonumber\\
=&v^a\Big|v_{ij}+d\frac{v_iv_j}{v}\Big|^2-\Big(d^2+\frac{q}{q-1}d+\frac{q}{2(q-1)}\Big)v^{a-2}|\nabla v|^4.
\end{align}
By means of the technique of Obata, we define a trace-free 2-tensor
\begin{align}\label{def E}
E_{ij}:=v_{ij}+d\frac{v_iv_j}{v}-\frac{1}{n}\Big(\frac{q}{q-1}+d\Big)\frac{|\nabla v|^2}{v}\delta_{ij}.
\end{align}
The tensor satisfies
\begin{align}\label{vijvij}
	\Big|v_{ij}+d\frac{v_iv_j}{v}\Big|^2=\Big|E_{ij}+\frac{1}{n}\Big(\frac{q}{q-1}+d\Big)\frac{|\nabla v|^2}{v}\delta_{ij}\Big|^2=|E|^2+\frac{1}{n}\Big(\frac{q}{q-1}+d\Big)^2\frac{|\nabla v|^4}{v^2}.
\end{align}
As what we have anticipated, \eqref{M-bN ob4} can be expressed as Obata type identities
\begin{prop}\label{Obata identity}
Let $v$ be a positive solution of \eqref{main equation}. The quantities $M,N$ are defined as \eqref{def of MN}. Then for $a=-\frac{q+1}{q-1}$ and $b,d$ satisfies $d=\frac{2q(1-b)-(q+1)}{2(q-1)}$, we have
\begin{align}\label{M-bN ob5}
	M-bN=v^a|E|^2+\bigg(\frac{1}{n}\Big(\frac{q}{q-1}+d\Big)^2-\Big(d^2+\frac{q}{q-1}d+\frac{q}{2(q-1)}\Big)\bigg)v^{a-2}|\nabla v|^4.
\end{align}
\end{prop}

\subsection{Choice of parameter b}
As what we have mentioned, \eqref{c geq 0} hopes $d\geq 0$. However, to make Obata's skill work, we apply
\begin{prop}\label{c leq 0}
For $n\geq 3$, $1<q<\frac{n}{n-2}$, we have
\begin{align}\label{ob geq 0}
\frac{1}{n}\Big(\frac{q}{q-1}+d\Big)^2-\Big(d^2+\frac{q}{q-1}d+\frac{q}{2(q-1)}\Big)\geq 0,
\end{align}
provided 
\begin{align}
-\frac{2q}{(q-1)(q+1)}\leq d\leq 0.
\end{align}
\end{prop}
\noindent{\it Proof.} Based on the observation that the left-hand side of \eqref{ob geq 0} is decreasing on $n$. It suffices to show the case when $n\to \frac{2q}{q-1}$, i.e.
\begin{align}
\frac{q-1}{2q}\Big(\frac{q}{q-1}+d\Big)^2-\Big(d^2+\frac{q}{q-1}d+\frac{q}{2(q-1)}\Big)\geq 0,\nonumber
\end{align}
equivalently,
\begin{align}
-\frac{q+1}{2q}d^2-\frac{1}{q-1}d\geq 0.\nonumber
\end{align}
we complete the proof.
\\\qed\\
\indent Combining with \eqref{c geq 0} and Proposition \ref{c leq 0}, it is sure that we need $d=0$, i.e. $b=\frac{q-1}{2q}$. In this case, we may both eliminate the term $\int_{S^{n-1}}v^{a-1}v_n^2$ in \eqref{Integration on Q-2} and make sure the coefficient of $v^{a-2}|\nabla v|^4$ in \eqref{M-bN ob5} is positive. To conclude what we have stated, from \eqref{M-bN ob5} and \eqref{Integration on Q-2},
\begin{cor}
Let $v$ be a positive solution of \eqref{main equation}. The quantities $M,N$ are defined as \eqref{def of MN}. Then for $a=-\frac{q+1}{q-1}$ and $b=\frac{q-1}{2q}$, we have
\begin{align}\label{M-bN ob}
	M-bN=v^a|E|^2+\frac{q}{2n(q-1)^2}\Big(n-(n-2)q\Big)v^{a-2}|\nabla v|^4\geq 0,
\end{align}
and 
\begin{align}\label{Integration on Q-4}
\begin{aligned}
	\int_{B^n}(M-bN)=&\frac{\lambda(3q-5)-2}{2}\int_{S^{n-1}}v^a|\nabla_{S^{n-1}}|^2\\
	&-(n-2)(q-2)\frac{q+1}{2}\lambda\int_{B^n}v^a|\nabla v|^2\\
	&+\bigg(\lambda \frac{q+1}{2}+\frac{(n-2)q-n}{2}\bigg)\int_{S^{n-1}}v^av_n^2.
\end{aligned}
\end{align}
\end{cor}
\begin{rmk}\label{rmk}
The proof of Proposition \ref{c leq 0} also suggests that for $a=-\frac{q+1}{q-1}$ and any $b'\in \bR$ satisfies
\begin{align}\label{range b'}
\frac{q-1}{2q}\leq b'\leq \frac{q^2+4q-1}{2q(q+1)},
\end{align}
we always have
\begin{align}
M-b'N\geq 0,
\end{align}
where $v$ is a positive solution of \eqref{main equation} and $M,N$ are defined as \eqref{def of MN}. In fact, if we define
\begin{align}
d':=\frac{2q(1-b')-(q+1)}{2(q-1)},
\end{align}
then correspondingly, the range of $d'$ is precisely
\begin{align}\label{range d'}
-\frac{2q}{(q-1)(q+1)}\leq d'\leq 0.
\end{align}
It is easy to check that \eqref{range d'} is equivalent to \eqref{range b'}.
\end{rmk}

\section{Auxiliary function and Improved identities}\label{sec:4}
To deal with more complicate situations, an introduction of auxiliary function $\phi$ here is essential. Set
\begin{align}\label{phi-1}
\phi(x)=\frac{|x|^2+c}{2}
\end{align}
with $c>0$ to be determined. It is evident to see
\begin{align}\label{phi-2}
\phi_i=x_i,\quad \phi_{ij}=\delta_{ij}\quad\text{and}\quad \phi(x)\in \Big[\frac{c}{2},\frac{1+c}{2}\Big) \quad\text{in }B^n,
\end{align}
and
\begin{align}\label{phi-3}
\phi(x)\equiv \frac{1+c}{2},\quad \text{and}\quad \nabla \phi=\nu,\quad \text{on }S^{n-1}. 
\end{align}

\subsection{Improved inequalities involving auxiliary function}
We define
\begin{align}
P:=\divg(v^a|\nabla v|^2\nabla \phi),\qquad Q:=\divg(v^a\la\nabla v,\nabla \phi\ra\nabla v),
\end{align}
where $a=-\frac{q+1}{q-1}$. Direct computations lead
\begin{align}
\label{def P}P=&2v^av_{ij}v_ix_j-\frac{q+1}{q-1}v^{a-1}|\nabla v|^2v_ix_i+nv^a|\nabla v|^2,\\
\label{Compute PQ} Q=&v^av_{ij}v_ix_j-\frac{1}{q-1}v^{a-1}|\nabla v|^2v_ix_i+v^a|\nabla v|^2.
\end{align}
If we focus on the quantity $(q-3)P+4Q$, then by \eqref{def P}-\eqref{Compute PQ},
\begin{align}\label{PQ identity}
(q-3)P+4Q=&2(q-1)v^av_{ij}v_ix_j-(q-1)v^{a-1}|\nabla v|^2v_ix_i+\Big(4-(3-q)n\Big)v^a|\nabla v|^2.
\end{align}
Recall the definition of $M,N$ in \eqref{def of MN}, and note clearly $\phi\equiv \frac{1+c}{2}$ on $S^{n-1}$. If we multiply $\phi$ on $M,N$ and integrate them on $B^n$ respectively, with the help of divergence theorem and \eqref{phi-2}-\eqref{phi-3}, we will obtain
\begin{align}\label{M identity}
\int_{B^n}M\phi=&\int_{B^n}(v^av_{ij}v_i)_j\phi=\int_{B^n}(v^av_{ij}v_i\phi)_j-\int_{B^n}v^av_{ij}v_ix_j\nonumber\\
=&\int_{S^{n-1}}v^av_{in}v_i\phi-\int_{B^n}v^av_{ij}v_ix_j\nonumber\\
=&\frac{1+c}{2}\int_{S^{n-1}}v^av_{in}v_i-\int_{B^n}v^av_{ij}v_ix_j\nonumber\\
=&\frac{1+c}{2}\int_{B^n}M-\int_{B^n}v^av_{ij}v_ix_j,
\end{align}
and
\begin{align}\label{N identity}
\int_{B^n}N\phi=&\int_{B^n}(v^a\Delta vv_i)_i\phi=\int_{B^n}(v^a\Delta vv_i\phi)_i-\int_{B^n}v^a\Delta vv_ix_i\nonumber\\
=&\int_{S^{n-1}}v^a\Delta vv_n\phi-\frac{q}{q-1}\int_{B^n}v^{a-1}|\nabla v|^2v_ix_i\nonumber\\
=&\frac{1+c}{2}\int_{S^{n-1}}v^a\Delta vv_n-\frac{q}{q-1}\int_{B^n}v^{a-1}|\nabla v|^2v_ix_i\nonumber\\
=&\frac{1+c}{2}\int_{B^n}N-\frac{q}{q-1}\int_{B^n}v^{a-1}|\nabla v|^2v_ix_i.
\end{align}
Putting \eqref{M identity}-\eqref{N identity} into the integration of \eqref{PQ identity} and noting $b=\frac{q-1}{2q}$, we arrive at a neat result that
\begin{align}\label{PQ identity-3}
\int_{B^n}\Big((q-3)P+4Q\Big)=& 2(q-1)\int_{B^n}(M-bN)\Big(\frac{1+c}{2}-\phi\Big)-\Big((3-q)n-4\Big)\int_{B^n}v^a|\nabla v|^2.
\end{align}

On the other hand, if we right use divergence theorem on the integration of $P,Q$ in $B^n$, we will derive by use of \eqref{phi-3} that
\begin{align}
\label{int M}\int_{B^n}P=&\int_{S^{n-1}}v^a|\nabla v|^2\la \nabla \phi,\nu\ra=\int_{S^{n-1}}v^a|\nabla v|^2=\int_{S^{n-1}}v^a|\nabla_{S^{n-1}}v|^2+\int_{S^{n-1}}v^av_n^2,\\
\label{int N}\int_{B^n}Q=&\int_{S^{n-1}}v^a\la \nabla v,\nabla \phi\ra\la\nabla v,\nu\ra=\int_{S^{n-1}}v^av_n^2.
\end{align}
Combine with \eqref{int M} and \eqref{int N}, we conclude
\begin{align}\label{PQ identity-4}
\begin{aligned}
\int_{B^n}\Big((q-3)P+4Q\Big)=&(q-3)\int_{S^{n-1}}v^a|\nabla_{S^{n-1}}v|^2+(q+1)\int_{S^{n-1}}v^av_n^2.
\end{aligned}
\end{align}
Thus \eqref{PQ identity-3} and \eqref{PQ identity-4} imply
\begin{align}\label{PQ identity-5}
\begin{aligned}
\int_{B^n}(M-bN)\phi=&\frac{1+c}{2}\int_{B^n}(M-bN)+\frac{3-q}{2(q-1)}\int_{S^{n-1}}v^a|\nabla_{S^{n-1}}v|^2\\
&-\frac{(3-q)n-4}{2(q-1)}\int_{B^n}v^a|\nabla v|^2-\frac{q+1}{2(q-1)}\int_{S^{n-1}}v^av_n^2.
\end{aligned}
\end{align}
Putting \eqref{Integration on Q-4} into \eqref{PQ identity-5}, we obtain an improved identity
\begin{cor}\label{main cor}
Let $v$ be a positive solution of \eqref{main equation}. The quantities $M,N$ are defined as \eqref{def of MN} where $a=-\frac{q+1}{q-1}$, $b=\frac{q-1}{2q}$. Set $\phi(x)=\frac{|x|^2+c}{2}$ with $c>0$ to be determined. Then we have
\begin{align}\label{final identity}
\begin{aligned}
\int_{B^n}(M-bN)\phi=&\bigg(\frac{1+c}{2}\Big(\frac{\lambda(3q-5)-2}{2}\Big)+\frac{3-q}{2(q-1)}\bigg)\int_{S^{n-1}}v^a|\nabla_{S^{n-1}}v|^2\\
&-\bigg(\frac{1+c}{2}(n-2)(q-2)\frac{q+1}{2}\lambda+\frac{(3-q)n-4}{2(q-1)}\bigg)\int_{B^n}v^a|\nabla v|^2\\
&+\bigg(\frac{1+c}{2}\Big(\lambda \frac{q+1}{2}+\frac{(n-2)q-n}{2}\Big)-\frac{q+1}{2(q-1)}\bigg)\int_{S^{n-1}}v^av_n^2.
\end{aligned}
\end{align}
\end{cor}

\section{Proofs of Theorem \ref{main thm}}\label{sec:5}
After having addressed all the preceding conclusions, we are now in a position to establish the proof of Theorem \ref{main thm}. Integrating \eqref{Compute PQ} in $B^n$ and combining with \eqref{int N}, \eqref{M identity} and \eqref{N identity}, it is obvious that
\begin{align}\label{ineq vn^2 vn}
\int_{S^{n-1}}v^av_n^2=\int_{B^n}Q=\int_{B^n}(M-\frac{1}{q}N)(\frac{1+c}{2}-\phi)+\int_{B^n}v^a|\nabla v|^2.
\end{align}
Choose $c=\frac{2}{\lambda(q-1)}-1\geq 1$ in Corollary \ref{main cor}. In this case, \eqref{final identity} becomes
\begin{align}\label{final identity2}
\begin{aligned}
\int_{B^n}(M-bN)\phi=&\bigg(1-\frac{1}{\lambda(q-1)}\bigg)\int_{S^{n-1}}v^a|\nabla_{S^{n-1}}v|^2\\
&-\frac{1}{2}\Big((n-2)q-n\Big)\int_{B^n}v^a|\nabla v|^2\\
&+\frac{(n-2)q-n}{2\lambda(q-1)}\int_{S^{n-1}}v^av_n^2.
\end{aligned}
\end{align}
Putting \eqref{ineq vn^2 vn} into \eqref{final identity2} to eliminate $\int_{B^n}v^a|\nabla v|^2$, we obtain
\begin{align}\label{final identity3}
\begin{aligned}
\int_{B^n}(M-bN)\phi=&\frac{1}{2}((n-2)q-n)\int_{B^n}(M-\frac{1}{q}N)(\frac{1+c}{2}-\phi)\\
&+\bigg(1-\frac{1}{\lambda(q-1)}\bigg)\int_{S^{n-1}}v^a|\nabla_{S^{n-1}}v|^2\\
&-\frac{(n-(n-2)q)(1-\lambda(q-1))}{2\lambda(q-1)}\int_{S^{n-1}}v^av_n^2.
\end{aligned}
\end{align}
Recall when $n\geq 3$, we have $1<q\leq \frac{n}{n-2}\leq 3$, Thus
\[
\frac{q-1}{2q}\leq \frac{1}{q}\leq \frac{q^2+4q-1}{2q(q+1)}.
\]
By Remark \ref{rmk}, we conclude that 
\[
M-\frac{1}{q}N\geq 0.
\]
So the conditions $1<q\leq \frac{n}{n-2}$ and $0<\lambda\leq \frac{1}{q-1}$ imply the right-hand side of \eqref{final identity3} is no greater than $0$. Combining with \eqref{M-bN ob} that the left-hand side of \eqref{final identity3} is no less than $0$, we conclude both sides of \eqref{final identity3} equal to $0$. Then the continuity of non-negative function $(M-bN)\phi$ implies
\begin{align}
(M-bN)\phi\equiv 0,\qquad \forall x\in B^n.
\end{align}
Taking into account that the range of $\phi$ strictly exceeds $0$, we attain
\begin{align}\label{M-bN==0}
M-bN\equiv 0,\qquad \forall x\in B^n.
\end{align}
Again by \eqref{M-bN ob}, if $1<q<\frac{n}{n-2}$, we derive $|\nabla v|\equiv 0$ in $B^n$. Thus $v$ is constant. \\
\indent On the other hand, if $0<\lambda<\frac{1}{q-1}$, then the right-hand side of \eqref{final identity3} equals to $0$ forces
\[
\int_{S^{n-1}}v^a|\nabla_{S^{n-1}}v|^2\equiv 0.
\]
Thus $v$ is constant on $S^{n-1}$. So $u$ is a harmonic function with constant boundary value, which must be constant by maximum principle.\\
\indent It remains to discuss the case $q=\frac{n}{n-2}$ and $\lambda=\frac{1}{q-1}$. In this case, \eqref{M-bN==0} and \eqref{M-bN ob} only implies $|E|=0$. Putting $a=-\frac{q+1}{q-1}$ and $b=\frac{q-1}{2q}$ into the definition \eqref{def E} of $E$, we have for any $1\leq i,j\leq n$,
\begin{align}\label{E=0}
v_{ij}=\frac{\Delta v}{n}\delta_{ij}.
\end{align}
By taking derivative on \eqref{E=0} and taking a sum, we derive for any $1\leq i\leq n$,
\[
(\Delta v)_i=\sum_{j=1}^nv_{jji}=\sum_{j=1}^nv_{ijj}=\sum_{j=1}^n\Big(\frac{\Delta v}{n}\delta_{ij}\Big)_j=\frac{1}{n}\sum_{j=1}^n(\Delta v)_j\delta_{ij}=\frac{1}{n}(\Delta v)_i.
\]
Since $n\geq 3$, we have $(\Delta v)_i=0$ for any $1\leq i\leq n$. Thus $\Delta v$ is constant. Suppose $\Delta v\equiv 2nr$ for some $r\in \bR$, then from \eqref{E=0}, we obtain
\begin{align}\label{expression of v}
v(x)=r|x|^2+\la \zeta,x\ra+s\qquad\forall x\in B^n,
\end{align}
for some $s\in \bR$ and $\zeta\in \bR^n$. Direct computations imply
\begin{align}\label{nabla v}
\nabla v=2rx+\zeta\qquad \forall x\in B^n.
\end{align}
Putting \eqref{expression of v} and \eqref{nabla v} into \eqref{main equation}, we get the relations between $r,s,\zeta$ that
\begin{align}
\label{rs1} 4rs=|\zeta|^2,\\
\label{rs2} r=s-\frac{2}{n-2}.
\end{align}
Since $v$ is a positive funciton, necessarily, we have $v(0)>0$, i.e. $s>0$. Combining with \eqref{rs1} and \eqref{rs2}, we solve
\begin{align}\label{r,s}
r=\sqrt{\Big(\frac{1}{n-2}\Big)^2+\frac{1}{4}|\zeta|^2}-\frac{1}{n-2},\qquad s=\sqrt{\Big(\frac{1}{n-2}\Big)^2+\frac{1}{4}|\zeta|^2}+\frac{1}{n-2}.
\end{align}
Using the fact that for any $\zeta\in \bR^n$, there is a unique $\xi\in B^n$, s.t.
\begin{align}\label{xi}
\zeta=-\frac{4}{n-2}\frac{\xi}{1-|\xi|^2}.
\end{align}
Putting \eqref{xi} into \eqref{r,s}, we derive
\begin{align}\label{r,s2}
r=\frac{2}{n-2}\frac{|\xi|^2}{1-|\xi|^2},\qquad s=\frac{2}{n-2}\frac{1}{1-|\xi|^2}.
\end{align}
Then putting \eqref{xi} and \eqref{r,s2} into \eqref{expression of v}, we obtain
\begin{align}
v(x)=\frac{2}{n-2}\frac{1+|\xi|^2|x|^2-2\la \xi,x\ra}{1-|\xi|^2},
\end{align}
i.e.
\begin{align}
u(x)=\bigg(\frac{n-2}{2}\frac{1-|\xi|^2}{1+|\xi|^2|x|^2-2\la \xi,x\ra}\bigg)^{\frac{n-2}{2}},
\end{align}
for some $\xi\in B^n$. We complete Theorem \ref{main thm}.
\\\qed\\

Combing Guo-Hang-Wang \cite{GHW21} in $n=2$ and Theorem \ref{main thm} in $n\geq 3$, the conjecture is true. And thus we obtain a new proof of Beckner's inequality \eqref{Beckner inequality} by integral method.
\begin{cor}
The following inequalities hold
\[\begin{aligned}
|S^{n-1}|^{\frac{q-1}{q+1}}\bigg(\int_{S^{n-1}}u^{q+1}\bigg)^{\frac{2}{q+1}}\leq (q-1)\int_{B^n}|\nabla u|^2+\int_{S^{n-1}}u^2,\quad \forall u\in C^{\infty}(\bar B^n),
\end{aligned}\]
provided $1<q<\infty$, if $n=2$, and $1<q\leq \frac{n}{n-2}$, if $n\geq 3$.
\end{cor}
\noindent{\it Proof.} We define the Sobolev quotient of a function $u\in H^1(B^n)-\{0\}$ as
\begin{align}\label{def of Q}
Q_{\lambda,q}(u):=\frac{\int_{B^n}|\nabla u|^2+\lambda\int_{S^{n-1}}u^2}{\Big(\int_{S^{n-1}}|u|^{q+1}\Big)^{\frac{2}{q+1}}}.
\end{align}
For the case $1<q<\frac{n}{n-2}$, the trace operator $H^1(B^n)\to L^{q+1}(S^{n-1})$ is compact. Thus the minimization problem
\begin{align}\label{minimization problem}
S_{\lambda,q}:=\inf_{u\in H^1(B^n)}Q_{\lambda,q}(u)
\end{align}
is achieved by smooth positive functions which satisfies \eqref{initial equation}. By use of the conjecture, such minimizer $u$ must be
\begin{align}\label{const u}
u(x)\equiv \lambda^{\frac{1}{q-1}},\qquad \forall x\in B^n.
\end{align}
Putting \eqref{const u} into \eqref{def of Q}, since the constant achieves $S_{\lambda,q}$, we have for any $u\in C^{\infty}(\bar B^n)$,
\begin{align}
\frac{\int_{B^n}|\nabla u|^2+\lambda\int_{S^{n-1}}u^2}{\Big(\int_{S^{n-1}}|u|^{q+1}\Big)^{\frac{2}{q+1}}}\geq \lambda|S^{n-1}|^{\frac{q-1}{q+1}},\qquad \forall 0<\lambda\leq \frac{1}{q-1}.
\end{align}
Let $\lambda=\frac{1}{q-1}$, we obtain \eqref{Beckner inequality} while $1<q<\frac{n}{n-2}$.\\
\indent For the critical case $q=\frac{n}{n-2}$, we follow the method of continuity of Aubin and Trudinger, see for example \cite{GT83}. It suffices to show the function $q\mapsto S_{\lambda,q}$ is continuous on the left at $q=\frac{n}{n-2}$. To prove, for any $\epsilon>0$, $\exists u_1\in H^1(B^n)$, s.t. 
\begin{align}\label{S>Q}
S_{\lambda,\frac{n}{n-2}}\geq Q_{\lambda,\frac{n}{n-2}}(u_1)-\frac{\epsilon}{2}.
\end{align}
Since the function $q\mapsto Q_{\lambda,q}(u_1)$ is continuous in $(1,\frac{n}{n-2}]$. $\exists \delta>0$, s.t. $\forall q'\in (\frac{n}{n-2}-\delta,\frac{n}{n-2}]$, we have
\begin{align}\label{Q>Q}
Q_{\lambda,\frac{n}{n-2}}(u_1)\geq Q_{\lambda,q'}(u_1)-\frac{\epsilon}{2}.
\end{align}
Combining with \eqref{S>Q} and \eqref{Q>Q}, we obtain
\begin{align}\label{upper continuous}
S_{\lambda,\frac{n}{n-2}}\geq Q_{\lambda,\frac{n}{n-2}}(u_1)-\frac{\epsilon}{2}\geq Q_{\lambda,q'}(u_1)-\epsilon\geq S_{\lambda,q'}-\epsilon.
\end{align}
On the other hand, assume $u_2\in H^1(B^n)$ satisfies
\begin{align}\label{S'>Q'}
S_{\lambda,q'}\geq Q_{\lambda,q'}(u_2)-\frac{\epsilon}{2}.
\end{align}
By H\"older inequality, we have
\begin{align}\label{Holder Q}
Q_{\lambda,q'}(u_2)=\frac{\int_{B^n}|\nabla u_2|^2+\lambda\int_{S^{n-1}}u_2^2}{\Big(\int_{S^{n-1}}|u_2|^{q'+1}\Big)^{\frac{2}{q'+1}}}\geq \frac{\int_{B^n}|\nabla u_2|^2+\lambda\int_{S^{n-1}}u_2^2}{\Big(\int_{S^{n-1}}|u_2|^{\frac{2(n-1)}{n-2}}\Big)^{\frac{n-2}{n-1}}}|S^{n-1}|^{\big(\frac{n-2}{n-1}-\frac{2}{q'+1}\big)}.
\end{align}
Since the function $q'\mapsto |S^{n-1}|^{\big(\frac{n-2}{n-1}-\frac{2}{q'+1}\big)}$ is continuous in $q'$, we may decrease the above $\delta>0$, s.t. $\forall q'\in (\frac{n}{n-2}-\delta,\frac{n}{n-2}]$, the following inequalities hold
\begin{align}\label{continuous q'}
\frac{\int_{B^n}|\nabla u_2|^2+\lambda\int_{S^{n-1}}u_2^2}{\Big(\int_{S^{n-1}}|u_2|^{\frac{2(n-1)}{n-2}}\Big)^{\frac{n-2}{n-1}}}|S^{n-1}|^{\big(\frac{n-2}{n-1}-\frac{2}{q'+1}\big)}\geq S_{\lambda,\frac{n}{n-2}}|S^{n-1}|^{\big(\frac{n-2}{n-1}-\frac{2}{q'+1}\big)}\geq S_{\lambda,\frac{n}{n-2}}-\frac{\epsilon}{2}.
\end{align}
Combine with \eqref{S'>Q'}, \eqref{Holder Q} and \eqref{continuous q'}, we obtain
\begin{align}\label{lower continuous}
S_{\lambda,q'}\geq S_{\lambda,\frac{n}{n-2}}-\epsilon.
\end{align}
By \eqref{upper continuous}, \eqref{lower continuous} and the arbitrariness of $\epsilon$, the function $q\mapsto S_{\lambda,q}$ is continuous on the left at $q=\frac{n}{n-2}$. Thus
\begin{align}
S_{\lambda,\frac{n}{n-2}}=\lim_{q\to (\frac{n}{n-2})^-}S_{\lambda,q}=\lim_{q\to (\frac{n}{n-2})^-}\lambda|S^{n-1}|^{\frac{q-1}{q+1}}=\lambda|S^{n-1}|^{\frac{1}{n-1}}.
\end{align}
We finish the new proof of Beckner's inequality. 
\\\qed\\

\subsection*{Acknowledgment}
The authors were partially supported by NSFC Grant No.11831005. We would like to express our thanks to Yao Wan for his helpful discussions and suggestions.

%-------------------------------------------------
\end{document}